\documentclass[reqno]{amsart}
\usepackage{amscd} 
\usepackage{amsfonts} 
\usepackage{amssymb} 
\usepackage{latexsym} 

\newcommand{\ncm}{\newcommand}

\def\ME{\mathcal{E}}

\newtheorem{theorem}{Theorem}[section]
\newtheorem{prop}[theorem]{Proposition}

\newtheorem{cor}[theorem]{Corollary}

\newtheorem{lem&def}[theorem]{Lemma \& Definition}
\newtheorem{definition}[theorem]{Definition}
\newtheorem{example}[theorem]{Example}

\ncm{\End}{\mbox{\rm End}\,}

\def\Hom{\mbox{\rm Hom}\,}
\def\Add{\mbox{\bf Add}\, }
\def\FGP{\mbox{\bf FGP}\, }

\def\id{\mbox{\rm id}}

\def\into{\hookrightarrow}
\def\to{\rightarrow}

\def\o{\otimes}    %tensor product 
\def\b{\, \Box}    %cotensor product
\def\bra{\langle}
\def\ket{\rangle}

\ncm{\rarr}[1]{\stackrel{#1}{\longrightarrow}}
\ncm{\larr}[1]{\stackrel{#1}{\longleftarrow}}
\def\cop{\Delta}

\def\eps{\varepsilon}

\def\du1{\hat 1}

\def\-2{_{(-2)}}
\def\-1{_{(-1)}}
\def\0{_{(0)}}
\def\1{_{(1)}}
\def\2{_{(2)}}
\def\3{_{(3)}}
\def\4{_{(4)}}

\def\|{\, | \,}

\def\du1{\hat 1}
\def\lact{\triangleright}
\def\ract{\triangleleft}
\def\op{^{\rm op}}

\begin{document}

\title[Anchor maps and stable modules]{Anchor maps and
stable modules in depth two}
\author{Lars Kadison}
\address{Department of Mathematics \\University of Pennsylvania \\ David Rittenhouse
Laboratory \\ 209 South 33rd Street \\ Philadelphia, PA 19104-6395} 
\email{lkadison@c2i.net} 
\date{}
\thanks{}

\begin{abstract} 
An algebra extension $A \| B$ is right depth two if its tensor-square $A\o_B A$ is in the Dress category
$\Add {}_AA_B$.  We consider necessary
conditions for  right, similarly left, D2 extensions in terms of partial $A$-invariance of two-sided ideals in $A$
contracted to the centralizer.  Finite dimensional
algebras extending central simple algebras are shown
to be depth two.  
 Following P.~Xu, left and right bialgebroids
 over a base algebra $R$ may be defined in terms of anchor maps, or representations on $R$. The anchor maps for
the bialgebroids 
$S = \End {}_BA_B$ and $T = \End {}_AA\o_BA_A$
over the centralizer $R = C_A(B)$ are the
modules ${}_SR$ and $R_T$ studied in \cite{LK2005, KK,
LK0505}, which provide information about the bialgebroids and the extension \cite{LK2005c}. 
The anchor maps for the Hopf algebroids in \cite{KR,
LK0805} reverse the order of right multiplication and action
by a Hopf algebra element, and lift to the
 isomorphism in \cite{OP}.   
We sketch a theory of stable $A$-modules and their endomorphism rings  and generalize the smash product
decomposition in \cite[Prop.\ 1.1]{LK2003} to any $A$-module.
We observe that Schneider's coGalois theory in
\cite{HS} provides examples of codepth two, such
as the quotient epimorphism of a finite dimensional 
normal Hopf  subalgebra. A
homomorphism of finite dimensional coalgebras is codepth two if and only if  
its dual homomorphism of algebras is depth two.      
\end{abstract} 
\maketitle

\section{Introduction and preliminaries}

Anchor maps for bialgebroids are defined algebraically by Ping Xu \cite{PX}
based on quantization of certain triangular Lie bialgebroids.  From
this point of view, a classical cocommutative bialgebroid such as the univ.\
env.\ algebra of a Lie algebroid has an anchor map extending the anchor map
of a Lie algebroid, which is a bundle map from a real vector bundle $E$ over a smooth
manifold $X$ into the tangent bundle $T(X)$ \cite{SW}. For example, $T(X)$ itself
has bialgebroid with total algebra $\mathcal{D}(X)$, the algebra of diffential operators on $X$,
base algebra $C^{\infty}(X)$ with anchor map $\mu: \mathcal{D}(X) \to \End C^{\infty}(X)$
the usual action of differential operators on smooth functions.  
The anchor map of a Lie algebroid $E$ pulls back from $T(X)$ a considerable amount of
differential geometry such as Lie bracket, connection and De Rham cohomology \cite{SW}.  
Analogously, the anchor map of a bialgeboid 
encodes in the unit module of the associated monoidal category
  information about the bialgebroid:
we will provide some evidence for this,
mentioned in the abstract and treated in some detail in section~3.

Depth two extensions arise very naturally from a subalgebra pair $B \subseteq A$ satisfying
a certain projectivity condition on the tensor-square.  For example, a Galois
extension $A \| B$ over a projective $R$-bialgebroid $H$ is depth
two since $A \o_B A \cong A \o_R H$ as $A$-$B$-bimodules, where $B$ and $R$ are commuting
subalgebras within $A$ \cite{LK2005b}.  This in fact characterizes depth two extension or at least
its endomorphism ring extension \cite{LK2006}.  A simpler example occurs if $B$ is
in the center of $A$, then $B \subseteq A$ is depth two if $A$ is finite projective
as a $B$-module.  If we think of a central simple algebra $B$ as a ``noncommutative
point'' we might expect that any finite dimensional algebra $A$ extending $B$ be
depth two.  We provide a rigourous proof of this and a somewhat more general fact
in section~2.  

For a quantum subalgebra such as sub-group algebra, sub-Hopf algebras, twisted and skew
variants of these, the notion of depth two is closely related to, perhaps characterizes,  the notion of normal subobject.
In the Clifford theory of decomposition of induced modules from a normal subgroup
$H \ract G$, 
certain modules are stable over the normal subgroup, i.e.,
 are isomorphic to all  conjugate modules.
The endomorphism ring of the induced module of a stable module has
essentially the structure of a twisted group algebra over the quotient algebra $G/H$
by a result of Conlon \cite[11C]{CR}. Schneider \cite{HJS} extends this and other classical results
by Clifford, Green and Blattner, unified within the induced representation theory of Hopf-Galois
extensions.  In section~4 we take the point of view that f.g.\ Hopf-Galois extensions (such as the finite strongly
group-graded algebras or graded Clifford systems)  
are depth two and certain constructions in \cite{HJS} such as stable modules
will sensibly generalize.  The endomorphism ring of an induced, restricted left $A$-module $M$
is for example a smash product of the depth two right bialgebroid $T$
acting on $\End {}_BM$ over the centralizer $R = A^B$ (Theorem~\ref{th-endo}).

Finally by way of introduction, the question of whether a depth two Hopf subalgebra is
normal in a finite dimensional Hopf algebra leads to a notion of codepth two
homomorphism of coalgebras $C \to D$ in the author's paper \cite{LK2007},
since Hopf algebra homomorphisms are also homomorphisms of coalgebras.
As we show in section~5, a 
homomorphism of finite dimensional coalgebras is codepth two if and only if  
its dual homomorphism of algebras is depth two. 
Schneider's coGalois theory in \cite{HS} provides an answer to a question of
whether $H \to H/HK^+$ is codepth two for normal Hopf subalgebras, the
answer being supplied in section~5.

\subsection{Preliminaries.} We let $\Add M$ denote the Dress category of a module $M_C$
over a ring $C$, consisting of all $C$-modules isomorphic
to direct summands of finite direct sums $M \oplus \cdots \oplus M$, and all module homomorphisms between these.  We let $\FGP C$ denote the
category of finitely generated, projective right $C$-modules and all module homomorphisms between these.  
Recall that $\Add M$ is equivalent to the category
$\FGP \End M_C$ via the functor $X_C \mapsto \Hom_C(M,X)$, 
where $\Hom_C(M,X)$ is a right module over
$E= \End M_C$ via ordinary composition.  
(Its inverse functor is given by $- \otimes_E M$
where ${}_EM_C$ is the natural bimodule.)  

Let $A\| B$ be a unital associative algebra extension, such as subring $B \subseteq A$  with $1_B = 1_A$
or  a unital ring homomorphism $B \to A$. Let $k$
denote the ground ring, a field in the later sections~4 and~5. Note
that the natural $B$-$A$-bimodule $A$ is
 in $\Add A\o_B A$, since
the multiplication mapping $A \o_B A \to A$ is a $B$-$A$-split
epi. The same is true of ${}_AA_B \in \Add A \o_B A$.   

The converse condition defines the notion of depth two.  
The extension $A \| B$ is depth two (D2) if $A \o_B A  \in \Add A$
as natural $A$-$B$-bimodules (right D2) and $B$-$A$-bimodules (left D2).

Note that $\End {}_BA_A \cong C_A(B)$, the centralizer
of the extension, which we denote by $R = C_A(B)$,
via $r \mapsto \lambda_r$, left multiplication of $A$ 
by $r \in R$.  
There is then a category equivalence   $\Add A \cong \FGP R$; in particular,
$\Hom ({}_BA_A, {}_BA \o_B A_A) := T$
is finitely generated projective as a right $R$-module.   Note that 
\begin{equation}
\label{eq: te}
 T \cong (A \o_B A)^B \cong \End {}_AA \o_B A_A
\end{equation}
via $f \mapsto f(1)$, which we take as an identification of $T$ with the $B$-commutator
of $A \o_B A$, and $t \mapsto (a \o_B a' \mapsto
at^1 \o_B t^2 a')$ respectively, where we denote $t = t^1 \o_B t^2 \in T$
(notationally suppressing any summation over simple
tensors).  The last isomorphism induces the ring structure
on $T = (A \o_B A)^B$ given by
\begin{equation}
\label{eq: tee}
 uv = v^1 u^1 \o_B u^2 v^2, \ \ \ 1_T = 1_A \o_B 1_A.
\end{equation}

Given $\Add M_C = \Add N_C$, it follows that 
that $\End M_C$ and $\End N_C$ are Morita equivalent
(Hirata, 1968).
In particular, $R$ and $\End {}_BA \o_B A_A$
are Morita equivalent. The inverse equivalence
of course comes from
$$ \Hom ({}_BA\o_B A_A, {}_BA_A) \cong \End {}_BA_B := S$$
via $f \mapsto f(- \o_B 1_A)$.  Thus, a left D2 extension $A \| B$ has f.g.\ projective $R$-module structures
 ${}_RS$ and $T_R$ on the rings $S$ and $T$.  Similarly, ${}_RT$
and $S_R$ are finite projective $R$-modules for a 
right D2 extension $A \| B$.   

In case $A \| B$ is additionally a Frobenius extension,
we have  the algebraic structure of the Jones tower
of a type $II_1$ subfactor:
$$ B \to A \into A_1 \into A_2$$
where $A_1 = \End A_B \cong A \o_B A$ and
$A_2 = \End (A_1)_A \cong \End A \o_B A_A$.
In the case of depth two, the relative commutators
$R = C_A(B)$ and $C_{A_2}(B) = \End {}_BA \o_BA_A$
are Morita equivalent with context bimodules
$C_{A_1}(B) \cong \End {}_BA_B$ and $C_{A_2}(A) \cong
\End {}_AA\o_BA_A$, i.e. $S$ and $T$.  Thus the notion of 
depth two algebra extension recovers classical
depth two for subfactors (\cite{KS, KK} for further
details).  Below we examine a pairing between $S$ and $T$
that becomes \cite[8.9]{KS} the Szyma\`nski nondegenerate  pairing of $C_{A_1}(B)$
and $C_{A_2}(A)$ in \cite[(14)]{KN}, which transfers the algebra
structure of one centralizer to a coalgebra structure on the other when $R$ is
trivial.    

The following coordinates for left and right
depth two extensions are useful for concrete computation. Given a left D2 extension $A \| B$,
we have a split epi $A^n \to A \o_B A$ and thus 
 a left D2 quasibase 
 $ \beta_i \in S, t_i \in T $ satisfying in $A \o_B A$: 
 \begin{equation}
\label{eq: ld2qb}
  x \o y = \sum_{i=1}^n t_i \beta_i(x)y.
\end{equation}
Assuming a left D2 quasibase for an extension is equivalent to our defining condition
above: define a split epi $\pi: A^n \to A \o_B A$
by $(a_1, \ldots, a_n) \mapsto \sum_{i=1}^n  t_ia_i$
with section $\sigma: x \o y \mapsto (\beta_1(x)y, \ldots, \beta_n(x)y)$. (The corresponding $B$-$A$-endomorphism
$\sigma \circ \pi: A^n \to A^n$ is given by the idempotent matrix
$\underline{\mbox{\bf r}} = (r_{ij})$ in $M_n(R)$, where
$r_{ij} = \sum_{(t_j)} \beta_i(t_j^1)t_j^2$.)

Similarly, given a right D2 extension $A \| B$,
there is a right D2 quasibase: $\exists \gamma_j \in S, u_j \in T$ such that  
\begin{equation}
\label{eq: rd2qb}
 x \o y = \sum_{j=1}^m x\gamma_j(y)u_j.
\end{equation}
We will fix this notation for left and right D2 quasibases throughout this paper.  
\subsection{Mapping viewpoint on $R$-valued pairings between $S$ and $T$}
Given a algebra extension $A \| B$, recall the pairings of $S = \End {}_BA_B$
and $T = (A \o_B A)^B$ given by
$\bra \alpha, t \ket = \alpha(t^1)t^2$
and $[\beta, u ] = u^1 \beta(u^2)$
both with values in $R = C_A(B)$.  If
$A \| B$ is D2, these are nondegenerate
pairings inducing ${}_RS \stackrel{\cong}{\longrightarrow} {}_R\Hom (T_R, R_R)$ and $S_R \stackrel{\cong}{\longrightarrow} \Hom ({}_RT, {}_RR)_R$, respectively \cite{KS}.
  
Now note $\End {}_BA_B \stackrel{\cong}{\longrightarrow}
\Hom ({}_AA \o_B A \o_B A_A, {}_AA_A)$ via
$\alpha \mapsto \mu^2 \circ (\id_A \o \alpha \o \id_A)$,
with inverse $F \mapsto F(1 \o - \o 1)$.  Suppose
$F$ and $\alpha \in S$ are images of one another under
these mappings; suppose $G \in \End {}_AA \o_BA_A$ and $t \in T$ are images of 
one another under the mappings in eq.~(\ref{eq: te}) above.
There are two obvious $B$-linear mappings $A \o_B A \stackrel{\longrightarrow}{\to}
A \o_B A \o_B A$ given by $x \o y \mapsto 1 \o x \o y$ or $x \o y \o 1$.      
Then the two pairings $\bra \alpha , t \ket$ and $[ \alpha , t ]$
are equal to the two images of $1_A \o_B 1_A$ under composition of the following  mappings,
\begin{equation}
1 \o 1 \in 
A \o_B A \stackrel{G}{\to} A \o_B A \stackrel{\longrightarrow}{\to} A \o_B A \o_B A  \stackrel{F}{\to} A
\end{equation}
since $ \Hom ({}_BA \o_B A_B, {}_BA_B) \longrightarrow R $
via $H \mapsto H(1 \o 1)$.  

%%%%%%%%%%%%%%%%%%%%%%%%%%%%%%%%%%%%%%%%%%%%%%%%%%%%%%%%%%%%%%%%%%%%%%%%%%%%%%%%%%%%%%%%%%%%%%%%%%%%%%%%%%%%%%%%%%%%%
\section{Left and right depth two}

In \cite{LK0505} the author observed that the centralizer of a depth two extension is a normal
subring in Rieffel's sense.  Related to this, we point out a necessary condition for
an algebra extension to be left D2, resp.\ right
D2.   

\begin{prop}
Let $A \| B$ be an algebra extension with centralizer
$R = C_A(B)$. 
If  $A \| B$ is left D2,
then for each two-sided ideal $I \subseteq A$,
the ideal contracted to $R$ satisfies left partial
$A$-invariance:
\begin{equation}
A(I \cap R) \subseteq (I \cap R)A
\end{equation}
If $A \| B$ is right D2, then for each two-sided ideal $I \subseteq A$,
the ideal contracted to $R$ satisfies right partial
$A$-invariance:
\begin{equation}
 (I \cap R)A  \subseteq A(I \cap R) 
\end{equation} 
\end{prop}
\begin{proof}
Given $s \in I \cap R$ and $a \in A$ and a left D2 quasibase, we note
from eq.~(\ref{eq: ld2qb}) that
$$ as = \sum_i t^1_i s t^2_i \beta_i(a) \in (I \cap R)A$$ since $t_i \in (A \o_B A)^B$. 

Similarly, from a right D2 quasibase we obtain
$$ sa = \sum_j \gamma_j(a)u_j^1 s u^2_j \in A (I \cap R),$$ whence the second of the two set inclusions.   
\end{proof}

The conditions of left and right partial $A$-invariance may be compared with invariance under the left and right
adjoint actions ($\mbox{\rm ad}_{\ell}$- and $\mbox{\rm ad}_r$-invariance) of a Hopf subalgebra:  the identities  $a\1 \o \tau(a\2)a\3 = a \o 1$ and 
$a\1 \tau(a\2) \o a\3 = 1 \o a$, for an element $a$ in a Hopf algebra with antipode $\tau$,
are comparable to the eqs.~(\ref{eq: ld2qb}) and~(\ref{eq: rd2qb}) 
(cf.\ \cite{LK0505, Mo}).
With this criteria for one-sided D2 extensions, it is often easy
to identify certain extensions as not being D2.

\begin{example}
\begin{rm}
Suppose $U$ and $V$ are rings, while ${}_UM_V$
and ${}_VN_U$ are bimodules.  Form a special case
of the generalized matrix ring $A = \left(
\begin{array}{cc}
U & M\\
N & V
\end{array}
\right)$, a fiber product of upper and lower matrix algebras 
with multiplication given by
$$\left(
\begin{array}{cc}
u & m\\
n & v
\end{array}
\right)  \left(
\begin{array}{cc}
u' & m'\\
n' & v'
\end{array}
\right) = \left(
\begin{array}{cc}
uu' & um' + mv' \\
nu' + vn' & vv'
\end{array}
\right).$$
Let the subring $B = \left(
\begin{array}{cc}
U & 0\\
0 & V
\end{array}
\right)$.
Then the centralizer $R = \left(
\begin{array}{cc}
Z(U) & 0\\
0 & Z(V)
\end{array}
\right)$ where $Z(U)$ and $Z(V)$ are centers of $U$
and $V$.  Consider the two-sided ideal $I = \left(
\begin{array}{cc}
U & M\\
N & 0
\end{array}
\right)$ in $A$.   We note that $R \cap I = \left(
\begin{array}{cc}
Z(U) & 0\\
0 & 0
\end{array}
\right)$, $A(R \cap I) = \left(
\begin{array}{cc}
U & 0\\
N & 0
\end{array}
\right)$, and $(R \cap I)A = \left(
\begin{array}{cc}
U & M\\
0 & 0
\end{array}
\right)$. If $M \neq 0$ and $N \neq 0$,
then $A(R \cap I) \not \subset (R \cap I)A$
and $(R\cap I)A \not \subset A(R \cap I)$,
so $A \| B$ is neither left nor right D2 (respectively).

If $N = 0$ and $M \neq 0$, $A$ is an upper triangular
matrix and $B$ the ``diagonal'' subring. Choose instead
the two-sided ideal $I = \left(
\begin{array}{cc}
U & M\\
0 & 0
\end{array}
\right)$, show that $A(R \cap I) \not \subset (R \cap I)A$, so that $A \| B$ is not right D2.
It is not left D2 by  using instead the ideal
$J = \left(
\begin{array}{cc}
0 & M\\
0 & V
\end{array}
\right)$.
\end{rm}
\end{example}

It is an open problem if there exists a left D2 algebra extension which is not right D2
(or the reverse if we pass to opposite algebras).  The test above for the one-sided D2 property
might be helpful in finding such an algebra extension, if a certain 
extension showed signs of being depth two, with
 centralizer  intermediate in size and
the over-algebra having a sufficiently rich ideal structure.

We now turn to an example of depth two extension.

\begin{theorem}
Suppose $B$ is an Azumaya $k$-algebra and $A$ is a fin.\ gen.\
projective $k$-algebra containing $B$ as a subalgebra.
Then $A \supseteq B$ is a depth two extension.
\end{theorem}
\begin{proof}
Since $B$ is Azumaya, it is known (e.g. \cite[p.\ 46]{LK1999})
that there are  Casimir elements
$e_i \in (B \o_k B)^B$ and  elements $b_i \in B$
such that 
\begin{equation}
1 \o_k 1 = \sum_{i=1}^N e_i b_i.
\end{equation}
Then for all $b \in B, a \in A$
$$ b \o_k a = \sum_i e_i bb_i a = \sum_i f_i(g_i(b \o_k a))$$
where $g_i: B \o_k A \to A$ is defined by $g_i(b \o_k a) = bb_i a$ and $f_i: A \to B \o_k A$ is defined by $f_i(x) = e_i x$ for each $i = 1,\ldots,N$.  
It follows that
$$ {}_BB\o_k A_A \oplus * \cong {}_BA_A^N. $$

Since $B$ is separable as a $k$-algebra, $B^e$ is
a semisimple extension of $k$. (E.g.
if $k$ is a field, $B$ is semisimple and so
is $B^e = B \o_k B\op$.) Then ${}_BA_B$
is $k$-relative projective, and projective since
$A$ is projective over $k$.  Then ${}_BA_B \oplus * \cong {}_BB \o_k B_B^M$. 
Tensoring by $- \o_B A_A$, we obtain
$$ {}_BA \o_B A_A \oplus * \cong {}_BB \o_k A_A^M.$$
Combining the two displayed isomorphisms,
we obtain $A \o_B A \oplus * \cong A^{NM}$
as the natural $B$-$A$-bimodules, whence $A \| B$
is left D2.  It is similarly argued that $A \| B$
is right D2.  
\end{proof}

If $a_k \in A$ and $G_k \in \Hom (A, B^e)$ 
($k = 1,\ldots,M$) denote
a finite projective basis for $A$ as a right $B^e$-module, the proof above converts to the left D2 quasibase for $A \| B$: ($x,y \in A$)  
\begin{equation}
x \o_B y = \sum_{i,k} e^1_i a_k \o_B e^2_i (b_i \cdot
G_k(x))y
\end{equation}
where $\beta_{ik}(x) = b_i \cdot G_k(x)$ are in fact  
$B$-valued endomorphisms in $\End {}_BA_B$ and $t_{ik} = e^1_i a_k \o_B e^2_i$
are in $R \o 1 \subseteq (A \o_B A)^B$.   

We will study elsewhere  the
more general setting of composite extensions $A \|B \| C$, where 
$A \| B$, $A \| C$ are D2 extensions
and $B \| C$ is H-separable, the
total rings underlying the bialgebroids
$T_{A|B} = (A \o_B A)^B$ and
$T_{A|C} = (A \o_C A)^C$ are Morita
equivalent in an interesting way.
The Morita context bimodules are
$P = (A \o_B A)^C$ and $Q = (A \o_C A)^B$, admitting a calculus extending
that in \cite{KS} and the centralizers
$A^B$ and $A^C$  
are functorial images of one another
under tensoring by $P$ and $Q$.

%%%%%%%%%%%%%%%%%%%%%%%%%%%%%%%%%%%%%%%%%%%%%%%%%%%%%%%%%%%%%%%%%%%%%%%%%%%%%%%%%%%%%%%%%%%%%%%%%
\section{Bialgebroids in terms of anchor maps}

A bialgebroid $H$ over a base ring $R$ is usually defined as an $R^e$-ring, 
$R$-coring with grouplike element $1_H$
and an augmentation ring $\eps:H_H \to R_H$ with multiplicative coproduct,
a definition that stays closest to the usual definition
of a Hopf algebra as an augmented algebra, coalgebra with homomorphic coproduct.  However, there
is a slightly different, equivalent way of defining a bialgebroid which comes from the theory of Lie
algebroids, their universal enveloping algebras (which are cocommutative Hopf algebroids)
and quantized variants of these. This is Ping Xu's definition  \cite{PX} of a left $R$-bialgebroid in terms of an anchor mapping, instead of a counit, an anchor map being
a representation of $H$ on $R$ yielding the unit module in the tensor category of $H$-modules. In this section we  give this definition
for right bialgebroids and show some useful aspects of the anchor mapping.    

The definition of a right bialgebroid  
 with total algebra $H$ and base algebra $R$ 
in terms of an anchor mapping is the following (cf.\ \cite{BM, Lu, PX} for corresponding
definition of left bialgebroid):
\begin{enumerate}
\item $H$ and $R$ are unital, associative $k$-algebras,
\item there are commuting algebra homomorphisms $R \stackrel{\sigma}{\longrightarrow} H \stackrel{\tau}{\longleftarrow} R^{\rm op}$
(called \textit{source} and \textit{target},
respectively) where for each $r,s \in R$
we have $\sigma(r) \tau(s) = \tau(s) \sigma(r)$,  
\item fix the $R$-$R$-bimodule ${}_RH_R$ given by
$$ r \cdot h \cdot s = h \sigma(s) \tau(r) \ \ \ (h \in H) $$ 
\item $H$ has $R$-$R$-bilinear \textit{comultiplication} $\cop: H \to H \o_R H$
given by notation $\cop(h) = h\1 \o h\2$
where $\cop(1_H) = 1_H \o_R 1_H$. 
 
\item multiplicativity of $\cop$ with technical pre-condition: for all $h, g \in H$, $r \in R$, 
\begin{equation}
 h\1  \o_R \tau(r) h\2 = \sigma(r) h\1 \o_R h\2  
\end{equation} 
\begin{equation}
\cop(h g) = \cop(h) \cop(g)
\end{equation}

\item an (anchor) map $\mu: H \to \End_k R$, an $R$-$R$-bimodule morphism w.r.t.\ the $R$-$R$-bimodule $\End R$ given by  $r \cdot f \cdot s = rf(-)s$, an anti-homomorphism and
right action of $H$ on  $R$ (whence we write $\mu(h)(r) := r \ract h$) satisfying:
\item $h\1 \sigma( r \ract h\2) = \sigma(r)h $
\item $h\2 \tau( r \ract h\1) = \tau(r)h $.
\end{enumerate}
The existence of an anchor map is equivalent to  the 
existence of a counit $\eps: H \to R$ (cf.\ \cite{BM},
$(H, \cop, \eps)$ becomes an $R$-coring).
For example, set $\eps(h) = 1_R \ract h$,
which is $R$-bilinear since $\mu$ is, and 
note that  $\eps(1_H) = 1_R$,
$\eps(h\1) \cdot h\2 = h\2 \tau(1 \ract h\1) = h$
 and
$$ \eps(\sigma(\eps(h))g) = 1 \ract \sigma (\eps(h)) g
= (1 \ract h) \ract g = \eps(hg). $$

Conversely, given the counit $\eps: H \to R$ the
anchor is computed from
\begin{equation}
\mu(h)(r) = \eps(\sigma(r)h) = \eps(\tau(r)h).
\end{equation}

In the tensor category of right $H$-modules, the anchor
mapping is the unit module structure on $R$, so
that for each $H$-module $V$, $R \o_R V \cong V$ as well as $V \o_R R \cong V$ as $H$-modules.
The anchor map may also be viewed as an arrow into the terminal object $(\End R)\op$ in the category of $R$-bialgebroids.  
 
\begin{example}
\begin{rm}
Let $A \| B$ be a right D2 extension. We recall
that $T = (A \o_B A)^B$ is a right bialgebroid over the centralizer
$R = C_A(B)$ (\cite[5.2]{KS}, two-sided depth two is not needed in the argument). The ring structure on $T$ is given in eq.~(\ref{eq: tee}). Note that $\sigma(r) = 1 \o_B r$
and $\tau(s) = s \o_B 1$ define homomorphism and
anti-homomorphism $R \to T$ that commute at all values
and induce from the right the $R$-$R$-bimodule
${}_RT_R$ given by $s \cdot t \cdot r = st^1 \o_B t^2 r$.  The comultiplication $\cop: T \to T \o_R T$ given
by $\cop(t) = \sum_j (t^1 \o_B \gamma_j(t^2)) \o_R
u_j$ and anchor mapping
\begin{equation}
\mu(t)(r) = r \ract t = t^1 r t^2
\end{equation}
are $R$-bilinear and satisfy
$$\sigma(r)t = t^1 \o_B rt^2 = (t^1 \o_B \gamma_j(t^2))(1 \o_B u_j^1 r u_j^2) = t\1 \sigma(r \ract t\2) $$
and similarly $\tau(r)t = t\2 \tau(r \ract t\1)$. 
The counit is given by $\eps(t) = \mu(t)(1) = t^1 t^2$.
We note that the representation $\mu$ is the module algebra $R_T$ and studied in \cite{KK,LK0505} as a generalized Miyashta-Ulbrich action.
\end{rm}
\end{example}
 
The corresponding anchor mapping based definition of
left bialgebroid is the opposite of the definition
above, given in detail in \cite{BM}.  Again we
are interested in the example coming from a depth
two extension $A \| B$.  We recall the left bialgebroid
structure on the endomorphism ring $S = \End {}_BA_B$
over $R$.  The source and target mapping $R \to S \leftarrow R\op$ are provided by the standard left
and right multiplication mappings $\lambda_r: x \mapsto rx$ and $\rho_s: x \mapsto xs$ for $r,s \in R$.  Of
course, $\lambda_r \circ \rho_s = \rho_s \circ \lambda_r$
and ${}_RS_R$ is given by $$r \cdot \alpha \cdot s
= \lambda_r \circ \rho_s \circ \alpha = r \alpha(-)s$$
for $\alpha \in S$. 
The comultiplication $\cop: S \to S \o_R S$ is given
by
\begin{equation}
\cop(\alpha) = \sum_j \gamma_j \o_R (\alpha \ract u_j)
\end{equation}
where $\alpha \ract t = t^1 \alpha(t^2 -)$ is an action
of $T$ on $S$ discussed in more detail in the next section.  The counit $\eps: S \to R$ given by
$\eps(\alpha) = \alpha(1)$ together with $\cop$
provides the $R$-$R$-bimodule $S$ an $R$-coring structure. Then the anchor map $\mu : S \to \End R$
is given
by
\begin{equation}
\mu(\alpha)(r) = \eps(\alpha \circ \lambda_r)(1) = \alpha(r)
\end{equation}
This gives $R$ the structure of a left $S$-module algebra.  The underlying module ${}_SR$ has been
studied  in \cite{LK2005, LK0505}. 

As shown in the next example, a comparison of anchor maps may lift to an isomorphism
between the bialgebroids they represent.  

\begin{example}
\begin{rm}
Suppose $H$ is a Hopf algebra with antipode $\tau$ 
and $A$ a left $H$-module algebra. Consider
two left $A$-bialgebroid structures on the total
space $A \o_k A \o_k H$ which have been studied
recently.  
In \cite{LK0805} the left bialgebroid $A^e \bowtie H$ is defined (by considering the $S$ construction for
a (special depth two) pseudo-Galois extension) with multiplication given by  
\begin{equation}
(a \o b \bowtie h)(c \o d \bowtie k) = a (h\1 \lact c) \o d(\tau(k\2) \lact b) \bowtie h\2 k\1,
\end{equation}
See the paper \cite{LK0805} for the details; here, we will only need
to know that the source $s_L: A \to A^e \bowtie H$
is given by $s_L(a) = a \o 1_A \bowtie 1_H$, and the counit  by 
\begin{equation}
\eps(a \o b \bowtie h) = a (h \lact b).
\end{equation}

We compute the anchor map $\mu: A^e \bowtie H \to \End A$:
$$
\mu(a \otimes b \bowtie h)(x) = \eps((a \otimes b \bowtie h)(x \o 1_A \bowtie 1_H)) $$
$$ = 
\eps(a(h\1 \lact x) \o b \bowtie h\2)
= a(h \lact xb)
$$
whence ($a,b \in A, h \in H$) 
\begin{equation}
\mu(a \otimes b \bowtie h) = \lambda_a \circ (h \lact \cdot) \circ \rho_b.
\end{equation}

In the papers \cite{CM, KR} the same data inputs 
into a left bialgebroid  $A \odot H \odot A$
where multiplication is given by
\begin{equation}
(a \odot h \odot b)(c \odot k \odot d) = a(h\1 \lact c) \odot h\2 k \odot (h\3 \lact d) b
\end{equation} 
source homomorphism by $\lambda(a) = a \odot 1_H \odot 1_A$
and counit by
\begin{equation}
\eps(a \odot h \odot b) =  \eps_H(h)a b
\end{equation}
where $\eps_H$ is of course the counit of the Hopf algeba
$H$.  
The anchor for this bialgebroid is then ($a,b,x \in A$,
$h \in H$): 
$$\mu(a \odot h \odot b)(x) = \eps((a \odot h \odot b)
(x \odot 1_H \odot 1_A)  = \eps(a(h\1 \lact x) \odot h\2 \odot
b ) =
a (h  \lact x)b$$
i.e.
\begin{equation}
\mu(a \odot h \odot b) = \lambda_a \circ \rho_b \circ
(h  \lact \cdot) 
\end{equation}
Observe that 
$$
\lambda_a \circ (h \lact \cdot ) \circ \rho_b = \lambda_a \circ \rho_{h\2 \lact b} \circ (h\1 \lact \cdot ) ,
$$
which lifts to an isomorphism of bialgebroids 
$A^e \bowtie H \cong A \odot H \odot A$ given by 
\begin{equation}
 a \o b \bowtie h \longmapsto a \odot h\1 \odot h\2 \lact b 
\end{equation}
given in \cite{OP} (with inverse,
$ a \odot h \odot b \mapsto a \o \tau(h\2) \lact b \bowtie h\1 $, a bialgebroid homomorphism commuting
with source, target and counit maps and is
an $A$-coring homomorphism \cite{BW}).
\end{rm}
\end{example} 

%%%%%%%%%%%%%%%%%%%%%%%%%%%%%%%%%%%%%%%%%%%%%%%%%%%%%%%%%%%%%%%%%%%%%%%%%%%%%%%%%%%%%%%%%%%%%%%%%%%%%%%%%%%%%%%

\section{On stable modules and their endomorphism rings}

Suppose ${}_AM$ is a left $A$-module.  Let $\ME$ denote
its endomorphism ring as a module restricted to a  $B$-module: $\ME = \End {}_BM$.
There is a right action of $T$ on $\mathcal{E}$ 
given by $f \ract t = t^1 f(t^2 -)$ for $f \in \mathcal{E}$. This is
a measuring action and $\ME$ is a right $T$-module
algebra (as defined in \cite{KS, BW}), since
$$ (f \ract t\1)\circ (g \ract t\2) = \sum_i t^1_i f(t^2_i\beta_i(t^1)g(t^2 -)) = (f\circ g) \ract t. $$
The subring of invariants in $\mathcal{E}$ is $\End {}_AM$ since $\End {}_AM \subseteq \ME^T$ is obvious,
and $\phi \in \ME^T$ satisfies for $m \in M, a \in A$:
$$ \phi(am) = \sum_j \gamma_j(a) (\phi \ract u_j)(m)
= \sum_j \gamma_j(a) \eps_T(u_j) \phi(m) = a  \phi(m).$$
The next theorem  shows that the endomorphism
ring of the induced module is the smash product
ring of the bialgebroid $T$ with the endomorphism ring $\ME$
(generalizing \cite[1.1, $M = A$]{LK2003}), the isomorphism 
$\Psi: T \ltimes \mathcal{E} \longrightarrow \End\,_AA \o_B M$  being given by 
\begin{equation}
\Psi(t \o f)( a \o m) = at^1 \o_B t^2 f(m).
\end{equation}

\begin{theorem}
\label{th-endo}
Let $M$ be a left $A$-module and $\ME = \End {}_BM$.  If $A \| B$ is left depth two,
then there is a ring isomorphism $\Psi: T \ltimes \mathcal{E} \stackrel{\cong}{\longrightarrow}
 \End {}_AA \o_B M$. 
\end{theorem}
\begin{proof}
We let $\mu_M: A \o_B M \to M$ denote the multiplication mapping
 defined by $a \o m \mapsto am$. Letting $F \in \End\,_AA \o_B M$, define $\Phi:
\End {}_A A\o_B M \to T \o_R \ME$ by $$\Phi(F) = \sum_j t_j \o_R \mu_M \circ (\beta_j \o_B \id_M)F(1_A \o -).$$
Note that $\Phi \circ \Psi = \id$, since for $t \o f \in T \o_R \mathcal{E}$, 
$$ \sum_j t_j \o_R \mu_M(\beta_j \o \id_M)(t^1 \o t^2f(-)) = \sum_j t_j\beta_j(t^1)t^2 \o f = t \o f.$$
Next, given $F \in \End\,_AA \o_B M$, let $F^1(m) \o F^2(m) := F(1 \o m)$
noting that $F(a \o m) = aF^1(m) \o F^2(m)$. Observe that 
$\Psi \Phi = \id$ since
$$\Psi \Phi(F)(a \o m) = \sum_j a t^1_j \o t_j^2\beta_j(F^1(m))F^2(m) = 
aF^1(m) \o F^2(m) = F(a \o m).$$
Thus $\Psi$ is bijective linear mapping.

Verify that $\Psi$ is a ring isomorphism, using $\cop(u) = \sum_j t_j \o_R
(\beta_j(u^1) \o_B u^2)$:
\begin{eqnarray*}
\Psi((t \ltimes f)(u \ltimes g))(a \o m) & = & \Psi(t {u}\1 \ltimes (f \ract {u}\2)g)(a \o m) \\
& = & \sum_j a t_j^1 t^1 \o t^2 t^2_j \beta_j(u^1)f(u^2g(m)) \\
& = & a u^1 t^1 \o t^2 f(u^2 g(m)) \\
& = & \Psi(t \ltimes f) \circ \Psi(u \ltimes g)(a \o m). \qed 
\end{eqnarray*}
\renewcommand{\qed}{}\end{proof}

With $M = A$, we obtain the isomorphism,
\begin{equation}
\label{eq: duality}
 T \ltimes \End {}_BA \cong
 \End {}_AA \o_B A
\end{equation}
Note that if $A \| B$ is D2, the module ${}_AA \o_B A$ is finite projective and a generator
since 
$\mu: A \o_B A \to A$ splits as a left $A$-module epi.  Then
$\End {}_AA \o_B A$ is Morita equivalent
to $A$. If $A \| B$ is a Frobenius 
extension, $\End {}_BA$ is a smash product of $A$ and $S$.  In this case, eq.~(\ref{eq: duality}) may be viewed as
a duality result for D2 Frobenius extensions (cf.\ \cite{LK2003}, \cite[ch.\ 9]{Mo}). 

The theory of stable $B$-modules has the intent to
generalize a smash product result like the one above to a certain extent.   We sketch the beginnings of such a project by  extending the definition of stable modules and certain theorems in
Schneider \cite{HJS} to the bialgebroid-Galois extensions.  (Recall from \cite{LK2005b} that
such Galois extensions are characterized by being
D2 and balanced.)    

Suppose
that $A \| B$ is a faithfully flat, balanced, depth two
extension of algebras over a field $k$. Let $R$ again be the centralizer $C_A(B)$,
let $T\op$ be the left bialgebroid $(A \o_B A)^B$ over
$R$ with the opposite multiplication of that in eq.~(\ref{eq: tee}) and identical to $T$ as $R$-corings, 
and $BR$ the smallest subalgebra  in $A$ containing
$B$ (or the image of $B$) and $R$. We recall
that $A$ is a left $T\op$-comodule algebra with
coinvariants equal to $B$ (by faithful flatness
of the natural module ${}_BA$) \cite{LK2005}; in particular,
$A$ is a left $T$-comodule.  Of course,
$T$ is a left $T$-comodule over itself
(see \cite{BW} for comodules over corings).   

\begin{definition}
Suppose $M$ is a left $BR$-module. We say that
$M$ is \textit{$A$-stable} if
\begin{equation}
A \o_B M \cong T \o_R M
\end{equation}
by a left $B$-linear and left $T$-colinear isomorphism.
\end{definition}

This definition is most useful when $B = BR$, e.g. a maximal
commutative subalgebra of $A$ or a trivial centralizer $k1_A$ (then our bialgebroids are bialgebras
and Galois extensions are Hopf-Galois extensions).  If $BR = A$
we are in the situation below, that all $A$-modules are  $A$-stable. Recall our notation $t_i \in T$, $\beta_i 
\in S$ for a left D2 quasibasis.  

\begin{prop}
Any left $A$-module $M$ is $A$-stable via the isomorphism 
\begin{equation}
\Psi:\ A \o_B M \stackrel{\cong}{\longrightarrow} T \o_R M, \ \
a \o_B m \longmapsto \sum_i t_i \o_R \beta_i(a)m
\end{equation}
\end{prop}
\begin{proof}
We note that $\Psi$ is left $B$-linear
since $\Psi(ba \o_B m)= \sum_i t_i \o_R \beta_i (ba)m
= b \Psi(a \o_B m)$ since $\beta_i \in \End {}_BA_B$.
Recall that $A$ is a left $T$-comodule via
$\rho_L(a) = \sum_i t_i \o_R \beta_i(a)$,
and $T$ has coproduct $\cop(t) = \sum_i t_i \o_R (\beta(t^1) \o_B t^2)$.  Then
we compute that $\Psi$ is left $T$-colinear:
$$ (T \o \Psi) \circ (\rho_L \o M)(a \o m) =
\sum_{i,k} t_k \o t_i \o \beta_i(\beta_k(a))m $$
on the one hand, and 
$$ (\cop \o M) \circ \Psi(a \o m) = \sum_j t_j \o_R (\beta_j(t_i^1) \o_B t_i^2) \o_R \beta_i(a)m $$
on the other hand, equal elements of 
\begin{equation}
T \o_R T\o_R M
\stackrel{\cong}{\longrightarrow} A \o_B A\o_B M,\ \ t \o u \o m \longmapsto t^1 \o_B t^2 u^1 \o_B u^2m 
\end{equation}
since both map into the element $a \o_B 1_A \o_B m$.

Finally $\Psi$ has inverse mapping defined by
\begin{equation}
\Psi^{-1}: T \o_R M \stackrel{\cong}{\longrightarrow}
A \o_B M, \ \ t \o_R m \longmapsto t^1 \o_B t^2 m
\end{equation}
  which follows from eq.~(\ref{eq: ld2qb}) or
\cite[2.2]{LK2006}.  
\end{proof} 

Let ${}_BM$ be any $B$-module, $\ME = \End {}_BM$
its endomorphism ring, $N = A \o_B M$ its
induced $A$-module and $E = (\End {}_AN)\op$ its endomorphism
ring.  We make note of the natural module
$N_E$.  The depth two structure on $A \| B$ imparts
on $N$ an obvious left $T$-comodule structure with coaction $\cop_N$ enjoying a \textit{Hopf
module} compatibility condition w.r.t.\ the left $A$-module structure, since $A$ is a $T\op$-comodule algebra): ($a \in A, n \in N = A \o_B M$)
\begin{equation}
\cop_N(an) = a\-1 n\-1 \o a\0 n\0 \in T \o_R N
\end{equation}
where the coaction on $N$ is given by $n = a \o_B m \mapsto a\-1 \o_R a\0 \o_B m$.

Referring to our smash product decomposition above, we see that the  proposition below is automatic if $M$
is an $A$-module.

\begin{prop}
Let ${}_BM$ be a stable module.  Then there is a left $T\op$-comodule algebra structure on $E$
such that $N$ is a Hopf module w.r.t. $T$ and $E$
and the algebra monomorphism $\ME \into {}^{\rm co\, T}E$, $f \mapsto \id_A \o f$ is bijective.
\end{prop}
\begin{proof}
Define $\cop_E: E \to T \o_R E$ via the canonical isomorphisms (using hom-tensor relation and $T_R$
finite projective) where $h := \Hom (M, \cop_N)$: 
$$ E \cong \Hom({}_BM, {}_BN) \stackrel{h}{\longrightarrow} \Hom ({}_BM, T \o_R N) 
\cong T \o_R E.
$$
Denoting $\cop_E(F) = F\-1 \o F\0$, we note that
\begin{equation}
\label{eq: cope}
F\-1 \o_R F\0(1_A \o_B m) = \cop_N(F(1 \o m))
\end{equation}
Then $N$ is a (left-right) $(E,T)$-Hopf module since $$\cop_N(F(a \o m)) = \cop_N(a F(1 \o m)) =
a\-1 F\-1 \o_R F\0(a\0 \o m). $$
It follows similarly that $$\cop_E(F \circ G) = 
G\-1 F\-1  \o_R F\0 \circ G\0 $$ and
$\cop_E(1_E) = 1_T \o 1_E$ since $\sum_i t_i\beta_i(1_A) = 1 \o_B 1 = 1_T$.
Whence $E$ is a left $T$-comodule algebra. 
The mapping $\ME \to E$ is monic since $A_B$ is
faithfully flat.   
Clearly endomorphisms of the form $\id_A \o_B f$
for $f \in \ME$ are coinvariants of $\cop_E$
by eq.~(\ref{eq: cope}) and  that 
$\cop_N$ is $\rho_L \o \id_M$. For the converse, 
we first note that 
\begin{equation}
M \cong {}^{\rm co\, T}(A \o_B M)\ \  \mbox{\rm via} \ m
\mapsto 1 \o m
\end{equation}
since $T\o_R A$ is a Galois coring with coinvariants $B$ (cf.\ \cite{LK2005} and \cite[28.19]{BW}).
If $F\-1 \o F\0 = 1_T \o F$, it follows
from the displayed mapping 
that $F = \id_A \o_B g$ for some $g \in \ME$.    
\end{proof} 

If $T$ has an antipode satisfying a few 
axioms (e.g. \cite{BB}), one may moreover
 show that having a unitary and left $T$-colinear
mapping $J: T \to E$ is equivalent to $M$ being
isomorphic to a direct summand of an $A$-stable module
(cf.\ \cite[3.3]{HJS}).

%%%%%%%%%%%%%%%%%%%%%%%%%%%%%%%%%%%%%%%%%%%%%%%%%%%%%%%%%%%%%%%%%%%%%%%%%%%%%%%%%%%%%%%%%%%%%%%

\section{Hopf subalgebras and codepth two}

Let $C$ and $D$ be two coalgebras over a field $k$.  
The author defined a notion of codepth two for
a coalgebra homomorphism $C \to D$  \cite{LK2007}, which is dual to the notion of depth two for an 
algebra homomorphism $B \to A$. In this section we recall the definition 
of codepth two and provide an example
coming from Schneider's coGalois theory \cite{HS}.   
Let $H$ be a finite dimensional Hopf algebra over $k$.  
 A  Hopf subalgebra $K$
of $H$ has coideal $K^+ = \ker \eps_K$ and induces the coalgebra epimorphism $H \to H/HK^+$
which we observe to be codepth two in this section.

Let $g: C \to D$ be a homomorphism of coalgebras over a field $k$.  Then $C$ has an induced $D$-$D$-bicomodule structure given
by left coaction $$ \rho^L: C \to D \o C,\ \rho^L(c) = c\-1 \o c\0 := g(c\1) \o c\2,$$
and by right coaction $$ \rho^R: C \to C \o D,\ \rho^R(c) = c\0 \o c\1 := c\1 \o g(c\2). $$ These two coactions commute
by coassociativity; we denote the resulting $D$-$D$-bicomodule structure 
on $C$ by ${}^DC^D$ later in this section.  In a similar way, any $C$-comodule becomes a $D$-comodule via the homomorphism $g$, 
the functor of corestriction \cite[11.9]{BW}.  Unadorned tensors between modules are over $k$, we use a generalized Sweedler notation,
the identity is sometimes denoted by its object,  
and basic terminology such as coalgebra homomorphism, comodule or bicomodule is defined in the standard way such as in \cite{BW}. 

Recall that the cotensor product $$ C \b_D C = \{  c \o c' \in C \o C \, | \,
c\1 \o g(c\2) \o c' =  c \o g({c'}\1) \o {c'}\2 \},$$
where we suppress a possible summation $c \o c' = \sum_i c_i \o {c'}_i$.  
For example, if $g = \eps: C \to K$ the counit on $C$, $C \b_D C= C \o C$.  

Recall that $C \Box_D C$ is a natural $C$-$C$-bicomodule via the coproduct $\cop$ on $C$
applied as $\cop \o C$ for the left coaction and $C \o \cop$ for the right coaction \cite[11.3]{BW}.  
Then $\underline{\cop}:  \, C \to C\b_D C$ induced by $\cop$ (where $\underline{\cop}(c) := c\1 \o c\2$)
 is a $C$-$C$-bicomodule monomorphism.  As $D$-$C$-bicomodule it is split
by $c \o c' \mapsto \eps(c) c'$,
and as a $C$-$D$-bicomodule $\underline{\cop}$ is split by $c \o c' \mapsto c \eps(c')$
for $c \o c' \in C \b_D C$. (Since $c\1 \o g(c\2) \o c' = c \o g({c'}\1) \o {c'}\2$,
it follows that $g(c) \o c' = \eps(c) g({c'}\1) \o {c'}\2$,
whence $c \o c' \mapsto \eps(c)c'$ is left $D$-colinear.) 
For example, if $D = C$ and $g = \id_C$, then $C \b_C C \cong C$, since
$\underline{\cop}$ is surjective.     

It follows that $C$ is in general isomorphic to a direct summand of $C \b_D C$
as $D$-$C$-bicomodules: $C \b_D C \cong C \oplus * $. 
Left codepth two coalgebra homomorphisms have the special complementary property: 

\begin{definition}{\cite[6.1]{LK2007}}
\begin{rm}
A coalgebra homomorphism $g: C \to D$ is  \textit{left codepth two (coD2)}
if for some positive integer $N$, we have $D$-$C$-bicomodule isomorphism 
\begin{equation}
C \b_D C \oplus * \cong C^N, 
\end{equation}
i.e., the cotensor product $C \b_D C$ is isomorphic to a direct summand of
a finite direct sum of $C$ with itself as $D$-$C$-bicomodules.
Right codepth two coalgebra homomorphisms are similarly defined.
\end{rm}
\end{definition}

Let $D^* \to C^*$ be the  algebra extension
$k$-dual to $g: C \to D$.  Various comodule
structures also pass to modules over the dual
algebras.  
Left coD2 quasibases are given for each
$c \o c' \in C \b_D C$ by
\begin{equation}
\label{eq: lcd2qb}
 c \o c' = \sum_{i=1}^N \eta_i(c \o {c'}\1) \alpha_i({c'}\2) \o {c'}\3
\end{equation}
where $\eta_i \in (C \b_D C)^{*\, D^*}$
and $\alpha_i \in \End {}^DC^D$ are called \textit{left coD2
quasibases} for the coalgebra homomorphism $g: C \to D$
\cite{LK2007}. The equation is analogous to the eq.~(\ref{eq: ld2qb}).
There is a right bialgebroid structure on $\End {}^DC^D$ over the centralizer $C^{* \, D^*}$ \cite{LK2007}.  

Schneider introduces the following set-up in \cite{HS}
for a Hopf algebra $H$ with bijective antipode $\tau$, 
which we call coGalois coextension because it is dual to Galois $H$-extensions. Let $C$ be  a right $H$-module coalgebra.  This means that in addition to being
a coalgebra and right $H$-module, it satisfies
the obvious compatibility conditions:
\begin{equation}
(ch)\1 \o (ch)\2 = c\1 h\1 \o c\2 h\2
\end{equation}
and $\eps_C(ch) = \eps_C(c) \eps_H(h)$
for all $c \in C, h \in H$. Then there is the canonical coalgebra epi $p: C \to \overline{C} =
 C/ CH^+$ where $H^+ = \ker \eps_H$, the elements of vanishing counit, and $CH^+$ a coideal of $C$. 
We define $C \to \overline{C}$ to be \textit{coGalois}
in case the mapping 
\begin{equation}
\label{eq: can}
\mbox{\rm can}: C \o H \longrightarrow C \b_{\overline{C}} C, \ \ \ c \o h \longmapsto c\1 \o c\2 h 
\end{equation}
 is bijective. Note that can does indeed have codomain in the cotensor product 
since $h - \eps_H(h)1_H \in H^+$. 

\begin{prop}
Suppose $H$ is a finite dimensional Hopf algebra.
If the coalgebra epimorphism $p: C \to \overline{C}$
defined above is coGalois, then
it is left and right codepth two.
\end{prop}
\begin{proof}
It is not difficult to check that can is a left $C$-colinear and right $\overline{C}$-colinear homomorphism with respect to obvious $C$-$\overline{C}$-bicomodule structures on $C \o H$
and $C \b_{\overline{C}} \, C$ (via $\cop_C$
and $(\id \o p)\cop_C$). Let $\dim H = n$.
Then $C \b_{\overline{C}} \, C \cong C^n$ as $C$-$\overline{C}$-bicomodules, whence $C \to \overline{C}$ is right coD2.  

Consider the variant coGalois mapping
\begin{equation}
\mbox{\rm can}': C \o H \longrightarrow  C \b_{\overline{C}} C, \ \ \ c \o h \longmapsto c\1 h \o c\2 
\end{equation} 
This is easily checked to be left $\overline{C}$-colinear and right $C$-colinear.
But $\mbox{\rm can}' = \mbox{\rm can} \circ \Phi$
via the bijective mapping $\Phi: C \o H \to C \o H$
and its inverse 
defined by 
\begin{equation}
\Phi(c \o h) = ch\1 \o \tau(h\2), \ \ \ \Phi^{-1}(c \o h) = ch\2 \o \tau^{-1}(h\1)
\end{equation}
Whence $\mbox{\rm can}'$ is bijective, so $C \b_{\overline{C}} \, C \cong C^n$ as $\overline{C}$-$C$-bicomodules, whence $C \to \overline{C}$ is left coD2.
\end{proof}  

\begin{cor}
Let $H$ be a finite dimensional Hopf algebra 
and $K$ a  Hopf subalgebra of $H$.
Then the canonical coalgebra epimorphism
$p: H \to H/HK^+$ is codepth two.
\end{cor}
\begin{proof}
Follows from the proposition if we let $C = H$
be the obvious underlying right $K$-module coalgebra
(where $K$ takes the place of $H$ in the proposition). We note that $\mbox{can}: H \o K \to
H \b_{H/HK^+} \, H$ is split monic via the retract
$H \o H \to H \o H$, $x \o y \mapsto x\1 \o \tau(x\2)y$,
restricted to the image of can.  Schneider \cite[Theorem II]{HS} shows that the mapping
can in eq.~(\ref{eq: can}) is bijective if
injective, and $C$ is a projective right $H$-module.  But $H$ is free as a natural $K$-module
by the Nichols-Zoeller theorem \cite{Mo}.  It follows
that can is bijective, so that $p$ is codepth two.
\end{proof}
If $K$ is normal in $H$, i.e., it is $\mbox{\rm ad}_{\ell}$-stable $a\1 K \tau(a\2) \subseteq K$
 and $\mbox{\rm ad}_r$-stable
$\tau(a\1)K a\2 \subseteq K$ for all $a \in H$,
then $HK^+ = K^+H$ \cite[3.4.4]{Mo} so
$H/HK^+$ is the quotient Hopf algebra $\overline{H}$.  
\begin{cor}
Let $H$ be any Hopf algebra with bijective antipode
and $K$ a normal finite dimensional Hopf subalgebra.
Then the canonical epi $p: H \to \overline{H}$
is codepth two.
\end{cor}
\begin{proof}
This follows from the proposition and the proof
of the previous corollary, except that we use a result
of Schneider's \cite{HJS93} to conclude
that $H$ is free over $K$.  
\end{proof}

\subsection{Duality between codepth two
and depth two}  Let $k$ be a field. Suppose all $k$-algebras
and $k$-coalgebras are finite dimensional in this subsection. In this case, there is
a duality $M \mapsto M^*$ of finite dimensional $C$-$D$-bicomodules with finite dimensional $A$-$B$-bimodules,
where $A = C^*$ is the dual algebra
of $C$ (with convolution multiplication)
and $B = D^*$. The bimodule structure
is given by $(a \cdot m^* \cdot b)(m) =
a(m\-1) m^*(m\0) b(m\1)$ in the obvious
notation.    
 We next show that
a morphism of coalgebras $g: C \to D$ is codepth two
if and only if its dual morphism $g^*: B \to A$ of algebras is depth two.  Moreover, the
bialgebroid of a codepth two extension defined in \cite[6.9]{LK2007}
is anti-isomorphic to the bialgebroid
$S$ of the depth two dual algebra extension. 

Let $C$ and $D$ be finite dimensional coalgebras,
and $A = C^*$ and $B = D^*$ be their
dual algebras.  Of course,  $g: C \to D$ is
a coalgebra homomorphism if and only if $g^*: B \to A$, defined by $g^*(d^*) = d^* \circ g$
where $d^* \in D^* = \mbox{\rm Hom}_k (D,k)$,
is an algebra homomorphism.  
Let $R = C_A(B)$, the centralizer of 
the $B$-$B$-bimodule induced on $A$ by $g^*: B \to A$: thus $R = \{ c^* \in A \,
| \, \forall \, d^* \in B, g^*(d^*)c^* = c^* g^*(d^*) \} $.    

\begin{theorem}
Let $g: C \to D$ be a homomorphism of coalgebras.  Then $g: C \to D$ is (left) coD2 if and only if the algebra
homomorphism $g^*: B \to A$ is  D2.
Moreover, the $R$-bialgebroids $\End {}^DC^D$ and $\End {}_BA_B$ are anti-isomorphic.  
\end{theorem}
\begin{proof}
($\Rightarrow$) We are given a split
epimorphism of natural $D$-$C$-bicomodules
$C^N \to C \b_D C$.  Applying the duality
mentioned above, we obtain a split
monomorphism of natural $D^*$-$C^*$-bimodules $(C \b_D C)^* \to (C^*)^N$.  

  We note the 
\cite[Lemma 3.5]{I} which holds if  $C$
alone is finite dimensional: as $A$-$A$-bimodules, there is an isomorphism
$\pi$,  
\begin{equation}
\label{eq: Iovanov}
A \o_B A \stackrel{\cong}{\longrightarrow}  (C \b_D C)^*
\end{equation}
where $\pi(a \o c^*)(c \o d) = a(c)c^*(d)$
for $a, c^* \in A$, $c\otimes d \in C \b_D C$.  

Then $A \o_B A \oplus * \cong A^N$ as
$B$-$A$-bimodules.  Whence $g^*: B \to A$
is left D2.  Similarly, we argue that
$C \to D$ is right coD2 $\Rightarrow$
$g^*: B \to A$ is right D2.  

($\Leftarrow$) We are given a split
epi $A^N \to A \o_B A$ of $B$-$A$-bimodules.  Note that $A^* = C$,
so by dualizing we have a split monic
$(A \o_B A)^* \to C^N$ of $D$-$C$-bicomodules.  By eq.~(\ref{eq: Iovanov}), $(A \o_B A)^* \cong C \b_D C$
as $C$-$C$-bicomodules. By corestriction
then $C \b_D C \oplus * \cong C^N$ as
$D$-$C$-bicomodules.  Thus
$g: C \to D$ is left coD2.
Similarly, we argue that $g: C \to D$ is
right coD2 if $g^*: B \to A$ is right D2.

The right bialgebroid $\End {}^DC^D$
over $R$ described in \cite[6.9]{LK2007}
is anti-isomorphic to the left bialgebroid
$\End {}_BA_B$ described in section~3 via
the mapping ($\alpha \in \End {}^DC^D,
c^* \in A$)
\begin{equation}
\End {}^DC^D \rightarrow \End {}_BA_B, 
\ \ \ \alpha \longmapsto \hat{\alpha}, \ \mbox{\rm where} \ 
\hat{\alpha}(c^*) = c^* \circ \alpha
\end{equation}
We leave it as an exercise to show that
this defines an anti-isomorphism of $R$-bialgebroids.  For example, the transform of the target map on $r \in R$ is $\lambda_r$, since for each $c \in C$
$\widehat{t_R(r)}(\eta)(c) = r(c\1) \eta(c\2) = \lambda_r(\eta)(c)$
by equation \cite[(23)]{LK2007}. 
The counit $\eps_S$ of the transform of $\alpha \in
\End {}^DC^D$ is evaluation at $1_{C^*} = 
\eps_C$, 
$ \eps_S(\hat{\alpha}) = \hat{\alpha}(1_A) = \eps_C \circ \alpha = \eps_E(\alpha)$
by the counit equation \cite[(25)]{LK2007}. 
Moreover, $
\cop_S(\hat{\alpha}) = \widehat{\alpha\1} \otimes \widehat{\alpha\2}
$
since by equations \cite[(7)]{LK0805}, \cite[(32)]{LK2007} ($\phi, \eta \in C^*, c \in C$)
$$ (\hat{\alpha}\1(\phi) * \hat{\alpha}\2(\eta))(c) = \hat{\alpha}(\phi * \eta)(c) =  (\phi * \eta)(\alpha(c)) = \phi(\alpha(c)\1) \eta(\alpha(c)\2)
$$
$$ = \phi(\alpha\1(c\1)) \eta(\alpha\2(c\2)) = (\widehat{\alpha\1}(\phi) * \widehat{\alpha\2}(\eta))(c)$$
where $*$ represents the convolution product on $C^*$.  
   \end{proof} 

%%%%%%%%%%%%%%%%%%%%%%%%%%%%%%%%%%%%%%%%%%%%%%%%%%%%%%%%%%%%%%%%%%%%%%%%%%%%%%%%%%%%%%%%%%%%%%%%%%%%

%%%%%%%%%%%%%%%%%%%%%%%%%%%%%%%%%%%%%%%%%%%%%%%%%%%%%%%%%%%%%%%%%%%%%%%%%%%%%%%%%%%%%%%%%%%%%%%%

\end{document}